\documentclass[12pt]{article}
\usepackage{a4wide}
\usepackage{amsmath}
\usepackage{amsthm}
\usepackage{amsfonts}
\usepackage{amssymb}
\usepackage{epsfig}
\usepackage{cite}
\usepackage{caption}
\usepackage{diagbox}
\usepackage{units}
\usepackage{float}
\usepackage{hyperref}
\newtheorem{theorem}{Theorem}
\newtheorem{lemma}[theorem]{Lemma}
\newtheorem{proposition}[theorem]{Proposition}
\newtheorem{conjecture}[theorem]{Conjecture}

\begin{document}

\title{Cyclic Coloring of Plane Graphs with Maximum Face Size 16 and 17}

\author{Zden\v ek Dvo\v r\'ak\thanks{Computer Science Institute of Charles University, Prague, Czech Republic. E-mail: {\tt rakdver@iuuk.mff.cuni.cz}. Supported by the Center of Excellence project P202/12/G061 of the Czech Science Foundation and by the project LL1201 (Complex Structures: Regularities in Combinatorics and Discrete Mathematics) of the Ministry of Education of Czech Republic.}\and
				Michael Hebdige\thanks{Department of Computer Science, University of Warwick, Coventry CV4 7AL, UK. E-mail: {\tt m.hebdige@warwick.ac.uk}. The work of this author was supported by the Undergraduate Research Support Scheme of the University of Warwick and by the Leverhulme Trust 2014 Philip Leverhulme Prize of the fourth author.} \and 
        Filip Hl\'asek\thanks{Computer Science Institute of Charles University, Prague, Czech Republic. E-mail: {\tt fhlasek@gmail.com}. The work of this author has received funding from the European Research Council under the European Union's Seventh Framework Programme (FP/2007-2013) / ERC Grant Agreement n. 616787.}\and
				Daniel Kr\'al'\thanks{Faculty of Informatics, Masaryk University, Botanick\'a 68A, 602 00 Brno, Czech Republic, and Mathematics Institute, DIMAP and Department of Computer Science, University of Warwick, Coventry CV4 7AL, UK. E-mail: {\tt dkral@fi.muni.cz}. The work of this author was supported by the Leverhulme Trust 2014 Philip Leverhulme Prize and by the MUNI Award in Science and Humanities of the Grant Agency of Masaryk university.}\and
				Jonathan A.~Noel\thanks{Department of Mathematics and Statistics, University of Victoria, Victoria, BC, Canada, V8P 5C2. E-mail: {\tt noelj@uvic.ca}. This work was completed while the fifth author was a DPhil student at the University of Oxford.}}

\date{}
\maketitle

\begin{abstract}

Plummer and Toft conjectured in 1987 that the vertices of every $3$-connected plane graph with maximum face size $\Delta^\star$ can be colored using at most $\Delta^\star+2$ colors in such a way that no face is incident with two vertices of the same color. The conjecture has been proven for $\Delta^\star=3$, $\Delta^\star=4$ and $\Delta^\star\ge 18$. We prove the conjecture for $\Delta^\star=16$ and $\Delta^\star=17$.

\end{abstract}

\section{Introduction}

Problems concerning planar graphs have always been among the most extensively studied topics in graph theory.
In this paper,
we study a generalization of proper coloring introduced by Ore and Plummer in 1969~\cite{bib-ore69+}:
a {\em cyclic coloring} of a plane graph is a vertex-coloring such that any two vertices
incident with the same face receive distinct colors.
The Cyclic Coloring Conjecture of Borodin~\cite{bib-borodin84} asserts that every plane graph with maximum face size $\Delta^\star$ has a cyclic coloring with at most $\lfloor 3\Delta^\star/2\rfloor$ colors.
There are many results on the Cyclic Coloring Conjecture and related problems.
We would like to particularly mention the Facial Coloring Conjecture of~\cite{bib-kral05+},
which implies the Cyclic Coloring Conjecture for odd values of $\Delta^\star$.
This conjecture,
which was addressed e.g.~in~\cite{bib-havet10+,bib-havet+,bib-kral05+,bib-kral07+},
asserts that every plane graph has an $\ell$-facial coloring with at most $3\ell+1$ colors,
i.e.~a vertex coloring such that any vertices joined by a facial walk of size at most $\ell$ receive different colors.
We refer the reader to an excellent survey~\cite{bib-borodin13} by Borodin
for further results related to the Cyclic Coloring Conjecture.

Despite a significant amount of interest (see e.g.~\cite{bib-borodin92,bib-borodin99+,bib-havet+,bib-sanders01+,bib-zlamalova}),
the Cyclic Coloring Conjecture has been proven only for three values of $\Delta^\star$: the case $\Delta^\star=3$,
which is equivalent to the Four Color Theorem proven in~\cite{bib-appel77,bib-appel77b} (a simplified proof was given in~\cite{bib-robertson99+}),
the case $\Delta^\star=4$ known as Borodin's Six Color Theorem~\cite{bib-borodin84,bib-borodin95},
and the recently proven case $\Delta^\star=6$~\cite{bib-hebdige+}.
Amini, Esperet and van den Heuvel~\cite{bib-amini08+},
building on the work in~\cite{bib-havet07+,bib-havet08+},
proved an asymptotic version of the Cyclic Coloring Conjecture: for every $\varepsilon>0$,
there exists $\Delta_0$ such that every plane graph with maximum face size $\Delta^\star\ge\Delta_0$ has a cyclic coloring with at most $\left(\frac{3}{2}+\varepsilon\right)\Delta^\star$ colors.

The graphs which witness that the bound in the Cyclic Coloring Conjecture is the best possible contain vertices of degree two;
in particular, they are not $3$-connected.
In 1987, Plummer and Toft~\cite{bib-plummer87+} conjectured the following:

\begin{conjecture}[Plummer and Toft~\cite{bib-plummer87+}]
\label{conj}
Every $3$-connected plane graph with maximum face size $\Delta^\star$ has a cyclic coloring with at most $\Delta^\star+2$ colors.
\end{conjecture}

This conjecture is the main subject of this paper.
We remark that Conjecture~\ref{conj} fails if the assumption on $3$-connectivity is replaced with the weaker assumption that the minimum degree is at least $3$~\cite{bib-plummer87+}.
It is tight in the case $\Delta^\star=4$, as the triangular prism graph (the Cartesian product of $K_3$ and $K_2$) is $3$-connected, has maximum face size $4$ and cannot be cyclically coloured with fewer then $6$ colours.
However, we do not know of any tight example with maximum face size greater than four. 
In fact, the upper bound stated in Conjecture~\ref{conj} is known to not be tight for $\Delta^\star$ sufficiently large: Borodin and Woodall~\cite{bib-borodin02+} showed that every $3$-connected plane graph with maximum face size $\Delta^\star\ge 122$ has a cyclic coloring with at most $\Delta^\star+1$ colors,
and this constraint on $\Delta^\star$ was lowered to $\Delta^\star\geq60$ by Enomoto,
Hor{\v n}{\'a}k and Jendrol'~\cite{bib-enomoto01+}. 

Conjecture~\ref{conj} has been proven for all but finitely many values of $\Delta^\star$.
The cases $\Delta^\star=3$ and $\Delta^\star=4$ follow from Four Color Theorem and Borodin's Six Color Theorem,
respectively.
The conjecture was proven for $\Delta^\star\ge 61$ in~\cite{bib-borodin02+},
for $\Delta^\star\ge 40$ in~\cite{bib-hornak00+},
for $\Delta^\star\ge 24$ in~\cite{bib-hornak99+} and finally for $\Delta^\star\ge 18$ in~\cite{bib-hornak10+}.
Our main result is a proof of the cases $\Delta^\star = 16$ and $\Delta^\star = 17$ of Conjecture~\ref{conj}:

\begin{theorem}
\label{main-thm}
Every $3$-connected plane graph with maximum face size $\Delta^\star\in \{16,17\}$ has a cyclic coloring that uses at most $\Delta^\star+2$ colors.
\end{theorem}

We employ the discharging method to prove Theorem~\ref{main-thm};
we refer to the surveys~\cite{bib-cranston13+,bib-cranston17+} for a detailed exposition of the method.
We start by identifying a set of configurations that cannot be contained in a minimal counterexample,
i.e., a counterexample with the smallest number of vertices, in Section~\ref{sec-redu}.
Such configurations are referred to as \emph{reducible configurations}.
We then consider a minimal counterexample $G$ and assign \emph{initial charges} to the vertices and faces of $G$
with the property that the sum of the initial charges is negative.
We then redistribute the charge using \emph{discharging rules},
which are described in Section~\ref{sec-rules}.
The redistribution preserves the overall sum of the charges.
Finally, we show that if $G$ contains none of the reducible configurations, then
every vertex and face has non-negative charge after applying the rules in Section~\ref{sec-analysis}, which is a contradiction.

Unfortunately, the arguments related to checking the reducibility of some of the configurations and
the analysis of the final charge turned out to be complex and
we had to resort to computer assistance.
We have made our programs verifying the correctness of our proof
available on-line at~\url{http://www.ucw.cz/~kral/cyclic-16/};
we have also uploaded their source codes to arXiv as ancillary files.

\section{Notation}

In this section, we briefly review the notation used in our proof.
Throughout this paper, all graphs are plane graphs unless explicitly stated.
A \emph{$k$-vertex} is a vertex of degree $k$.
We also define an \emph{$(\geq k)$-vertex} to be a vertex with degree at least $k$, and
an \emph{$(\leq k)$-vertex} to be a vertex with degree at most $k$.
The \emph{size} of a face $f$ of a plane graph, denoted by $|f|$,
is the number of vertices that are incident with it.
Analogous to the definition of a $k$-vertex, a \emph{$k$-face} is a face of size $k$.
Similarly, a \emph{$(\geq k)$-face} and a \emph{$(\leq k)$-face} are faces that
have size at least $k$ and at most $k$, respectively.
The \emph{boundary walk} of a face in a plane graph is a sequence of vertices that bounds the face. A pair of distinct vertices are 
said to be \emph{cyclically adjacent} if they are incident with a common face. 
The \emph{cyclic degree} of a vertex is the number of vertices which are cyclically adjacent to it.

\begin{figure}
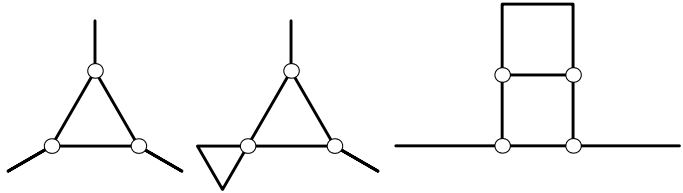

\begin{center}
\epsfbox{cyclic-16.38} \epsfbox{cyclic-16.39} \epsfbox{cyclic-16.40}
\end{center}
\caption{An A-triangle, a B-triangle and a column.}
\label{fig-triangles+column}
\end{figure}

Most of the configurations are depicted in Figures \ref{fig-triangles+column}--\ref{fig-gen} using the notation that we now describe.
A circled vertex in a configuration depicts its exact degree,
i.e., the vertex must be incident with as many edges as depicted in the figure.
When we need to name some of these vertices, we write their identifiers inside the circles.
The vertices depicted by bold circles
are required to have the cyclic degree equal to the number given in the figure next to the vertex
in addition to having the degree as depicted by the number of incident edges.
In addition, we sometimes restrict the face sizes by writing the constraint on the face size in the middle of the face;
for example in the first configuration depicted in Figure~\ref{fig-gen},
the bottom middle face is required to have size at least $9$ and
the top left face is required to have size $\Delta^\star+6-\ell$,
where $\ell$ is the size of the bottom middle face.

When describing the discharging rules,
we use the following terminology, some of which is illustrated in Figure~\ref{fig-triangles+column}.
Let $v_1v_2$ be a part of the boundary walk of a face $f$.
With respect to a face $f$,
a triangle $T=v_1v_2v_3$ is an \emph{A-triangle} if $\deg(v_1)=\deg(v_2)=\deg(v_3)=3$,
a \emph{B-triangle} if $\deg(v_1)=\deg(v_2)=3$, $\deg(v_3)=4$, and
the neighbors $x_1$ and $x_2$ of $v_3$ distinct from $v_1$ and $v_2$ are adjacent, and
a \emph{C-triangle} if $T$ is neither an A-triangle nor a B-triangle.
If $v_1v_2$ is incident with a $4$-face $Q=v_1v_2v_3v_4$, $\deg(v_1)=\deg(v_2)=\deg(v_3)=\deg(v_4)=3$ and
$v_3v_4$ is incident with another $4$-face, we say that $Q$ is a \emph{column} (with respect to $f$).
If $v_1vv_2$ is a part of the boundary walk of a face $f$ and neither $v_1v$ nor $v_2v$ is contained in a $(\le\!4)$-face,
we say that $v$ is \emph{isolated} (with respect to $f$).
If $\deg(v)=4$ and $v$ is contained in a $(\le\!4)$-face $f'$ that does not share an edge with $f$,
then we say that $f'$ is the \emph{sink} of $v$;
otherwise, $v$ is the sink of itself.

\section{Construction of discharging rules}

A large part of our proof is computer-assisted.
However, the proof itself was also constructed in a computer-assisted way.
Once the types of the discharging rules are fixed, e.g.,
we have decided to transfer some amount of charge from a face of some given size $\ell$ to an incident A-triangle,
the conditions that the final charge of every vertex and face is non-negative become linear constraints.
More precisely,
there is a single variable for each rule type, and
a single linear constraint for each possible neighborhood structure of a vertex or a face that
is not excluded by reducible configurations.
Hence,
the amounts of charge transferred by individual rules can be determined by solving a linear program (or
it is determined that no such amounts for the given set of rule types using existing reducible configurations).
This natural approach has been used by various researchers earlier (however,
we were not able to locate the reference for its first use).

Let us give more details about how we have proceeded in the case of our proof.
First, it is not clear what types of discharging rules should be considered.
We started with rule types close to those in Subsections~\ref{sec-large}, \ref{sec-med} and \ref{sec-heavy} and
later added further rule types in an ad hoc way.
We then repeated the following steps.
We ran a linear program solver to determine a minimal set of infeasible constraints.
Each such constraint corresponds to a particular neighborhood structure (configuration).
To proceed further, it was necessary to either find a new reducible configuration that would exclude one of these configurations or
add a new rule type that would move charge inside the configuration.
After this, we reran the linear program solver to determine a new minimal set of infeasible constraints.
When the solver produced a solution, we found a possible set of discharging rules, i.e., a proof.

Since each rule type adds a new variable to the linear program,
it is necessary to be careful with adding new rule types to keep the linear program of manageable size.
For example, it would have been ideal to have rules of the types as those in Subsection~\ref{sec-addit} for all face sizes
but this would have resulted in a linear program too large to be solved in a reasonable amount of time.
As a compromise, we have started with rougher rules from Subsection~\ref{sec-med} and
combine them with finer rules from Subsection~\ref{sec-addit}.
Another concern might be that most linear program solvers (we have used the Gurobi solver) work in floating arithmetic;
however, the solution output by the solver can be rounded to rational values and checked with exact arithmetic computations.

\section{Reducible configurations}\label{sec-redu}

In this section, we describe reducible configurations that are used in our proof.
The reducible configurations in Figures \ref{triangle0}--\ref{fig-gen} are named using the following convention:
the name of the configuration refers to the size of a face that the configuration primarily concerns and
the subscript is used to distinguish different types of configurations related to faces of the same size.
In addition to the configurations presented in the figures,
there are two additional reducible configurations:
DEG is the configuration comprised of a single vertex with cyclic degree at most $\Delta^\star+1$, and
TFEDGE is the configuration comprised of a $3$-face and a $(\le\!4)$-face sharing an edge.
These two configurations are reducible by \cite[Lemma 3.1(e)]{bib-hornak99+} and \cite[Lemma 3.6]{bib-hornak10+}, respectively.
We will also need the following proposition
to justify the reducibility of some of our configurations.

\begin{proposition}[Halin~\cite{bib-halin69}]
\label{prop-3ver}
If $G$ is a $3$-connected graph with at least five vertices,
then every vertex of degree three is incident with an edge $e$ such that
the contraction of $e$ yields a $3$-connected graph.
\end{proposition}



\subsection{Configurations TRIANGLE}

\begin{figure}
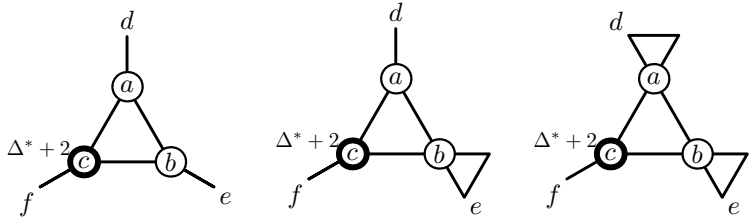

\begin{center}
\epsfbox{cyclic-16.47}
\hskip 4mm
\epsfbox{cyclic-16.48}
\hskip 4mm
\epsfbox{cyclic-16.49}
\end{center}
\caption{The configurations TRIANGLE$_0$.}
\label{triangle0}
\end{figure}

\begin{figure}
\begin{center}
\epsfbox{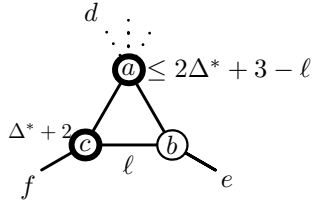}
\end{center}
\caption{The configuration TRIANGLE$_1$.}
\label{triangle1}
\end{figure}

\begin{figure}
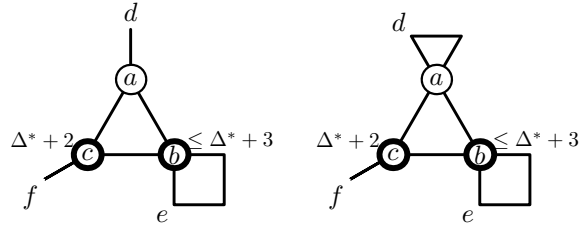

\begin{center}
\epsfbox{cyclic-16.51}
\hskip 4mm
\epsfbox{cyclic-16.52}
\end{center}
\caption{The configurations TRIANGLE$_2$.}
\label{triangle2}
\end{figure}

In this subsection, we introduce six reducible configurations referred to as TRIANGLE$_0$, TRIANGLE$_1$ and TRIANGLE$_2$, and
argue that they are indeed reducible.
The configurations can be found in Figures~\ref{triangle0}, \ref{triangle1} and \ref{triangle2}.
When analyzing the configurations, we use the following notation.
Let $abc$ be the $3$-face and let $A$ be the face that is incident with the vertices $a$ and $c$ but is not the $3$-face.
Let $A_{{v_1},\ldots,{v_n}}$ denote the set of colors that appear on all of the vertices that
are incident with $A$ except for the vertices $v_1,\ldots,v_n$.
Similarly, we define the faces $B$ and $C$ to be the faces incident with the edges $ab$ and $bc$, respectively.

\begin{lemma}\label{triangle-Lemma}
The three configurations denoted by TRIANGLE$_0$, which are depicted in Figure~\ref{triangle0}, and
the configuration denoted by TRIANGLE$_1$, which is depicted in Figure~\ref{triangle1}, and
the two configurations denoted by TRIANGLE$_2$, which are depicted in Figure~\ref{triangle2},
are reducible.
\end{lemma}

\begin{proof}
Assume that a minimal counterexample $G$ contains one of the configurations in Figures~\ref{triangle0},~\ref{triangle1} or~\ref{triangle2} and let $\ell=|C|$.
Note that $\ell\ge 5$ (otherwise, $G$ would contain the reducible configuration TFEDGE).
First contract the edge $cf$.
If the resulting graph were not $3$-connected,
then there would exist a vertex $x$ such that $c$, $f$ and $x$ form a vertex cut. 
Since $c$ is a $3$-vertex, the vertices $f$ and $x$ would form a vertex cut of size two in $G$, which is impossible.
Hence, the resulting graph is $3$-connected and
the minimality of $G$ implies that
the resulting graph has a cyclic coloring with $\Delta^\star+2$ colors.
This yields a coloring of $G$ of all vertices except for $c$.
If $A_c\cap C_{cf}$ is non-empty, then $c$ is cyclically adjacent to vertices of at most $\Delta^\star+1$ colors and
we can complete the coloring. Hence, assume that $A_c\cap C_{cf}=\emptyset$.
Without loss of generality,
we can assume that $a$ was colored with $1$, $b$ with $2$, $e$ with $3$, $f$ with $\ell$ and $d$ with $\Delta^\star+2$.
We can also assume that $C_{bcef}$ contains the colors from $4$ to $\ell-1$ and $A_{acdf}$ contains the colors from $\ell+1$ to $\Delta^\star+1$.

We first analyze the three configurations depicted in Figure~\ref{triangle0}.
If we can recolor $a$ with a color from $3,\dots,\ell-1$, then we can color $c$ with $1$, so $\{3,\dots, \ell-1\} \subseteq B_{abd}$.
If we can recolor $b$ with a color from $\ell+1,\dots,\Delta^\star+1$, then $c$ can be colored with $2$, and hence $\{\ell+1,\dots,\Delta^\star+1\} \subseteq B_{abd}$.
Therefore, $B_{abd}$ contains all the colors in $\{3,\dots,\Delta^\star+1\}\setminus\{\ell\}$,
which is impossible since $|B|\le\Delta^\star$.

We next analyze the configuration depicted in Figure~\ref{triangle1}.
If we can recolor $b$ with a color from $\ell+1,\dots,\Delta^\star+2$, then
we can color $c$ with $2$, and hence $\{\ell+1,\dots,\Delta^\star+2\} \subseteq B_{abe}$.
Likewise, if we can recolor $a$ with a color from $4,\dots,\ell-1$, then
we can color $c$ with $1$.
In particular, the vertex $a$ is cyclically adjacent to vertices with the colors from $4$ to $\ell-1$ and
to two vertices of each of the colors from $\{\ell+1,\dots,\Delta^\star+2\}$ (once on the face $A$ and once on the face $B$).
In addition to these $\ell-4+2(\Delta^\star-\ell+2)=2\Delta^\star-\ell$ vertices, $a$ is also cyclically adjacent to $b$, $c$, $e$ and $f$.
Hence, its cyclic degree must be at least $2\Delta^\star-\ell+4$,
which violates the description of the configuration.

It remains to analyze the configurations depicted in Figure~\ref{triangle2}.
Let $D$ be the set containing all colors in $B_{ab}$ and
the color assigned to the vertex of the $4$-face containing $b$ that is not adjacent to $b$;
we can assume that this color is $4$ without loss of generality.
Note that $b$ is the unique vertex contained in both $B$ and $C$ as, if it were not, 
then $G$ would fail to be $3$-connected.
Thus, since the cyclic degree of $b$ is at most $\Delta^\star+3$, it holds $|D|\le \Delta^\star+3-\ell$.
If we can recolor $b$ with a color from $A_{acf}$, then $c$ can be colored with $2$.
Hence, $A_{acf}\subseteq D$.
If $a$ can be recolored with $3$ or $4$, we can color $c$ with $1$.
Since this is impossible and $4\not\in A_c$, it holds that $\{3,4\}\subseteq D$.
We conclude that $D$ contains at least $|A_{acf}|+2=\Delta^\star+4-\ell$ colors,
which exceeds its size.
\end{proof}

\subsection{Computer assisted cases}

\begin{figure}
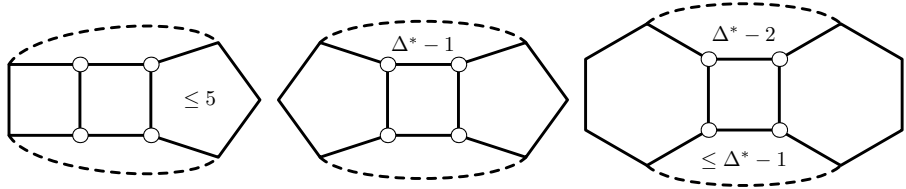

\begin{center}
\epsfbox{cyclic-16.14} \epsfbox{cyclic-16.15} \epsfbox{cyclic-16.16}
\end{center}
\caption{The configurations FOUR$_0$, FOUR$_1$ and FOUR$_2$.}
\label{fig-four}
\end{figure}

\begin{figure}
\begin{center}
\epsfbox{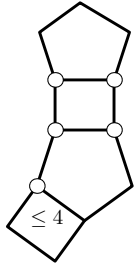}
\end{center}
\caption{The configuration FIVE.}
\label{fig-five}
\end{figure}

\begin{figure}
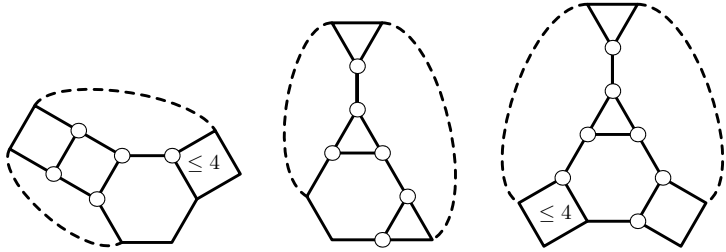

\begin{center}
\epsfbox{cyclic-16.8} \hskip 4mm \epsfbox{cyclic-16.43} \hskip 4mm \epsfbox{cyclic-16.44}
\end{center}
\caption{The configurations SIX$_1$, SIX$_2$ and SIX$_3$.}
\label{fig-six123}
\end{figure}

\begin{figure}
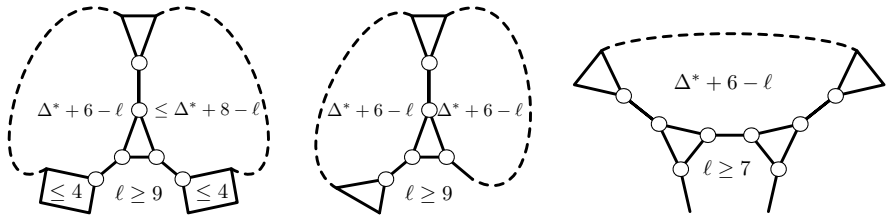

\begin{center}
\epsfbox{cyclic-16.28} \hskip 4mm \epsfbox{cyclic-16.29} \hskip 4mm \epsfbox{cyclic-16.30}
\end{center}
\caption{The configurations GEN$_0$, GEN$_1$ and GEN$_2$.}
\label{fig-gen}
\end{figure}

The remaining reducible configurations used in the proof are depicted in Figures~\ref{fig-four}--\ref{fig-gen}.
The configuration FOUR$_0$, which is depicted in Figure~\ref{fig-four},
is reducible by~\cite[Lemma 3.1(c)]{bib-hornak99+} and~\cite[Lemma 3.1(d)]{bib-hornak99+}.
The reducibility of the remaining configurations was verified with the assistance of a computer.
We have independently prepared two programs,
which are available at \url{http://www.ucw.cz/~kral/cyclic-16/} as
\texttt{test-reducibility1.c} and \texttt{test-reducibility2.cc}.
The input files needed to check the reducibility of the configurations are also available on-line.
Some of the configurations are checked using two input files,
corresponding to different settings described by the configuration.
For example,
there are two input files for the configuration FOUR$_1$ from Figure~\ref{fig-four},
one describing the setting where the bottom face is $(\le\!\Delta^\star-1)$-face and
the other the setting where it is a $\Delta^\star$-face.
Similarly, there are two input files for the configuration FIVE from Figure~\ref{fig-five},
one describing the setting where the bottom face is $3$-face and the other where it is a $4$-face.

\begin{figure}
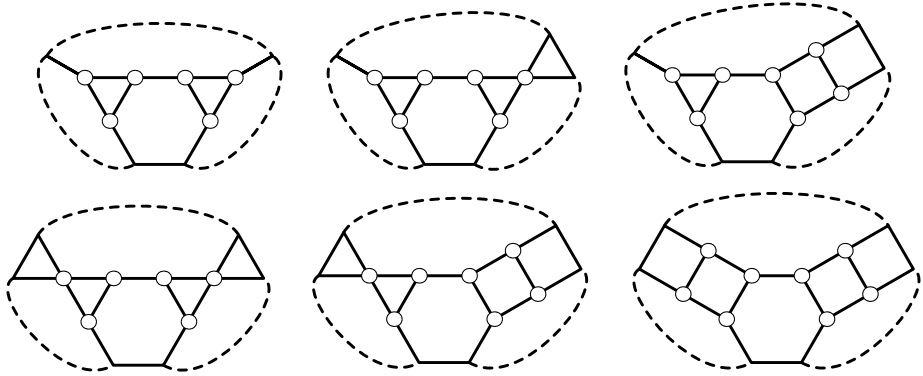

\begin{center}
\epsfbox{cyclic-16.1} \hskip 4mm \epsfbox{cyclic-16.2} \hskip 4mm \epsfbox{cyclic-16.4}
\vskip 2mm
\epsfbox{cyclic-16.3} \hskip 4mm \epsfbox{cyclic-16.5} \hskip 4mm \epsfbox{cyclic-16.6}
\end{center}
\caption{The configurations SIX$_0$.}
\end{figure}

\begin{figure}
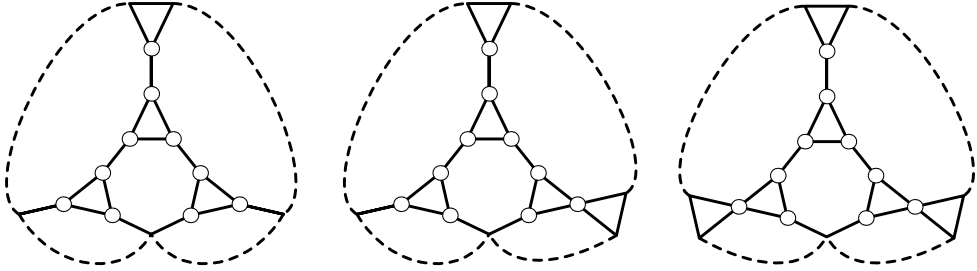

\begin{center}
\epsfbox{cyclic-16.7} \hskip 4mm \epsfbox{cyclic-16.9} \hskip 4mm \epsfbox{cyclic-16.10}
\end{center}
\caption{The configurations SEVEN$_0$.}
\end{figure}

\begin{figure}
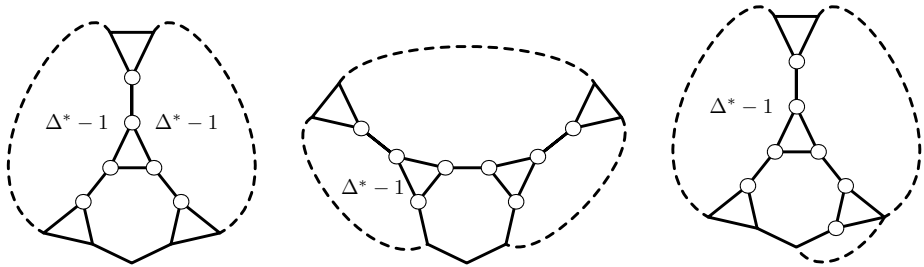

\begin{center}
\epsfbox{cyclic-16.11} \hskip 4mm \epsfbox{cyclic-16.12} \hskip 4mm \epsfbox{cyclic-16.45}
\end{center}
\caption{The configurations SEVEN$_1$, SEVEN$_2$ and SEVEN$_3$.}
\end{figure}

\begin{figure}
\begin{center}
\epsfbox{cyclic-16.13} \hskip 4mm \epsfbox{cyclic-16.17} 
\vskip 2mm
\epsfbox{cyclic-16.18} \hskip 4mm \epsfbox{cyclic-16.19}
\end{center}
\caption{The configurations EIGHT$_0$.}
\end{figure}

\begin{figure}
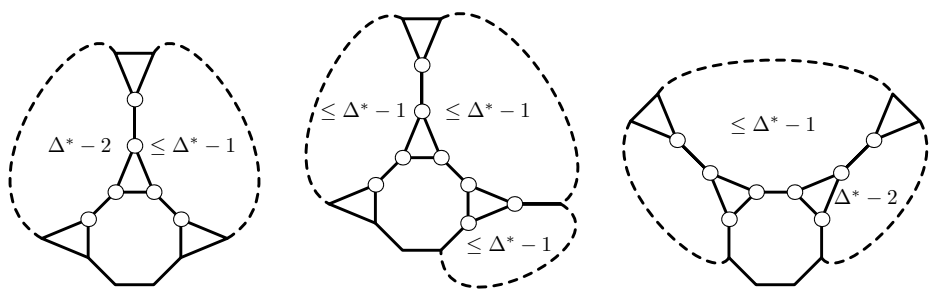

\begin{center}
\epsfbox{cyclic-16.20} \hskip 4mm \epsfbox{cyclic-16.21} \hskip 4mm \epsfbox{cyclic-16.22}
\end{center}
\caption{The configurations EIGHT$_1$, EIGHT$_2$ and EIGHT$_3$.}
\end{figure}

\begin{figure}
\begin{center}
\epsfbox{cyclic-16.23} \hskip 4mm \epsfbox{cyclic-16.24}
\end{center}
\caption{The configurations EIGHT$_4$.}
\end{figure}

\begin{figure}
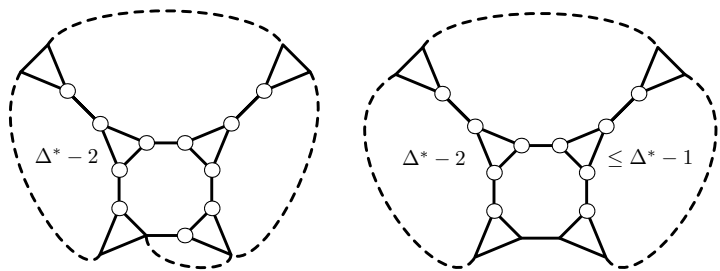

\begin{center}
\epsfbox{cyclic-16.25} \hskip 4mm \epsfbox{cyclic-16.26}
\end{center}
\caption{The configurations EIGHT$_5$ and EIGHT$_6$.}
\end{figure}

\begin{figure}
\begin{center}
\epsfbox{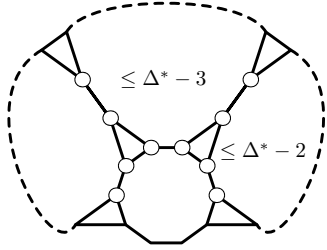} 
\end{center}
\caption{The configuration TEN.}
\end{figure}

We next describe the structure of the input files and the way used to reduce the configurations.
All configurations depicted in Figures~\ref{fig-four}--\ref{fig-gen} are verified to be reducible in the following manner. 
If a possible minimal counterexample $G$ contains the configuration in question,
we replace it with another configuration to obtain a graph $G'$ with a smaller number of vertices;
it will always be the case that the new configuration is obtained by contracting an edge and
this edge is drawn bold in Figures~\ref{fig-four}--\ref{fig-gen}.
Each input file consists of two blocks: the first block describes the new configuration and
the second block the original configuration, i.e., the configuration that we are verifying to be reducible.
The two blocks have a similar structure.
The first line of each block contains two integers $m$ and $n$.
The integer $m$ is the number of faces forming the configuration and
$n$ is the number of vertices with no neighbors outside the configuration (these are the circled vertices in the figures).
Let us call such vertices \emph{internal}.

\begin{figure}
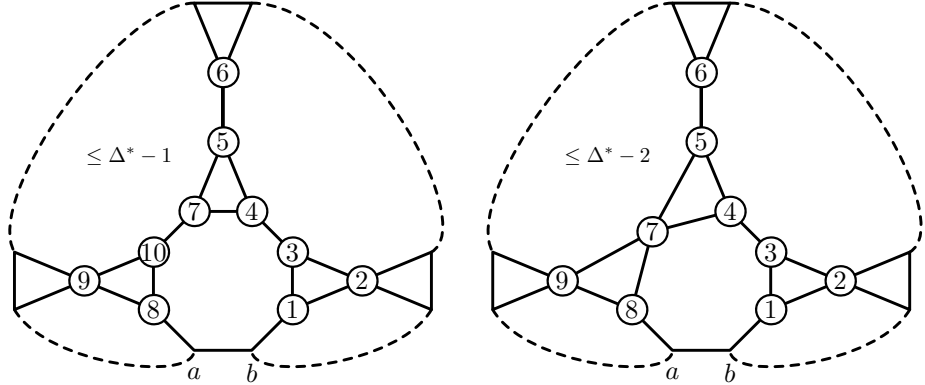

\begin{center}
\epsfbox{cyclic-16.53} \hskip 4mm \epsfbox{cyclic-16.56}
\end{center}
\caption{An example of how one of the configurations EIGHT$_0$ is reduced; the new configuration is on the right.}
\label{fig-explred}
\end{figure}

Each of the following $m$ lines describes one of the faces of the configuration.
There are two kinds of faces: \emph{bounded} faces with a specific size and
\emph{unbounded} faces with size between $\Delta^\star-c_1$ and $\Delta^\star-c_2$ (inclusively) for some $c_1$ and $c_2$.
A line describing a bounded face starts with $0$ and it is followed by the list of vertices incident with the face.
The vertices incident with the face that are internal with respect to the configuration
are represented by numbers between $1$ and $n$ and
the remaining vertices incident with the face are represented by lowercase letters;
the internal vertices have also their numbers shown in Figures~\ref{fig-four}--\ref{fig-gen}.
A description of an unbounded face in the first block starts with a range $a_1$--$a_2$; it is possible that $a_1=a_2$.
The rest of the line contains all internal vertices of the face and
possibly some additional vertices represented by lowercase letters (if there are no such additional vertices,
the line contains \verb|-|).
In addition to these vertices, the face is incident with $k$ vertices
where $k$ satisfies that $\Delta^\star+2-a_2\le k\le \Delta^\star+2-a_1$ (note that $a_i\not=c_i$ in general).
The values of $a_1$ and $a_2$ are determined using the constraint on the minimum cyclic degree (the configuration DEG),
the constraint on the maximum face size, and the configuration TRIANGLE$_0$.
For example, for the configuration SIX$_3$ depicted in Figure~\ref{fig-six123},
the two neighbors of the vertex $5$ are denoted by $e$ and $g$ in both input files for this configuration.
So, the bottom face is incident with at least $\Delta^\star-5$ additional vertices (otherwise,
the cyclic degree of the vertex $5$ would be at most $\Delta^\star+1$) and
at most $\Delta^\star-3$ additional vertices (otherwise, the bottom face would have size larger than $\Delta^*$);
so, the values of $a_1$ and $a_2$ are $5$ and $7$, respectively, and
the corresponding line of the input file is the following.
\begin{verbatim}
5-7 eg 5
\end{verbatim}
Similarly, since the vertex $3$ in Figure~\ref{fig-explred} has cyclic degree at least $\Delta^\star+3$,
the unbounded face incident with it must have size between $\Delta^\star-2$ and $\Delta^\star$.
Consequently, the number of non-internal vertices incident with this face is between $\Delta^\star-7$ and $\Delta^\star-5$,
which corresponds to the range $7$--$9$ given on the $4$-th line of the input file below.
In the second block,
the line describing an unbounded face starts with a positive integer giving the index of the corresponding face
in the first block (the indices start from one).
For example, the input file to verify the reducibility of the configuration EIGHT$_0$,
which is depicted in Figure~\ref{fig-explred},
is the following.

\begin{verbatim}
5 9
5-7 a 8,9
8-9 - 5,6,7,9
7-9 - 2,3,4,5,6
5-7 b 1,2
0 ab 1,3,4,7,8
5 10
1 a 8,9
2 - 5,6,7,9,10
3 - 2,3,4,5,6
4 b 1,2
0 ab 1,3,4,7,8,10
\end{verbatim}

\noindent Note that we have not specified the three $3$-faces formed by internal vertices,
e.g., the one formed by the vertices $1$, $2$ and $3$,
since the constraints that they impose on the coloring are implied by the presence of the other faces.
In addition, we have also not specified the existence of the other $3$-face containing the vertex $2$.
If the configuration can be checked to be reducible without this additional assumption,
it is also reducible with this additional assumption (we give more details further).

The program assumes the existence of a cyclic coloring of $G'$ using at most $\Delta^\star+2$ colors and
checks using this assumption that $G$ also has a cyclic coloring using at most $\Delta^\star+2$ colors.
When doing so, we assume that all the faces described in the input are pairwise different.
For example,
we assume that the face incident with $8$ and $9$ is different from the face incident with $1$ and $2$ in Figure~\ref{fig-explred}.
Each of the configurations that we analyze has a face $f$ such that
any other face described by the configuration shares a vertex with $f$.
Hence, all the faces described in the input can be assumed to be different because the graph $G$ is $3$-connected.

Another fact that needs to be verified is that the graph $G'$ is $3$-connected.
For most of our reductions, this is implied by Proposition~\ref{prop-3ver} in the following:
the graph $G'$ is obtained by contracting the edge incident with a $3$-vertex $v$ such that
the vertex $v$ is contained in a $3$-face formed by three $3$-vertices
but the contracted edge is not contained in the $3$-face (in the example,
the contracted edge joins the vertices $7$ and $10$).
Proposition~\ref{prop-3ver} guarantees that the $3$-vertex $v$ is incident with an edge that
can be contracted to obtain a $3$-connected graph, however,
contracting any of the other two edges incident with $v$ results in $2$-vertex.
In the remaining few cases, the $3$-connectivity of $G'$ follows by an easy analysis of the configurations.

We now describe how the program checks the existence of a cyclic coloring of $G$.
The program enumerates all possible colorings of non-internal vertices and
checks whether the coloring extends in $G'$, and if so,
it also checks that it extends in $G$.
Note that some of the colorings of non-internal vertices considered by the program are not feasible.
For example, we have neglected in the considered configuration one of the $3$-faces containing the vertex $2$ and
the constraints that it imposes.
Since testing the extendibility of a larger set of colorings does not harm the validity of our arguments,
this does not affect the correctness of our arguments as long as all the constraints on the coloring of internal vertices are represented.
In fact, this negligence is useful in the considered case
since the very same input file can be used to justify the reducibility of all the configurations EIGHT$_0$.

\section{Discharging rules}
\label{sec-rules}

In this section, we describe the discharging phase of our proof.
Each vertex $v$ of a minimal counterexample is assigned charge $\deg(v)-4$ and each face $f$ is assigned $|f|-4$.
Using Euler's formula, the overall sum of the initial charge is $-8$.
The charge then gets redistributed among the vertices and faces as follows.

First,
each $3$-vertex $v$ that is contained in exactly one $(\le\!4)$-face gets $1$ from this face, and
each $3$-vertex $v$ that is contained in two $4$-faces $f_1$ and $f_2$ gets $\frac{1}{2}$ from each of $f_1$ and $f_2$;
note that a $3$-vertex cannot be contained in a $3$-face and a $(\le\!4)$-face.
Other rules are more complex and are described in the rest of the section.
We start with simpler rules to redistribute the charge, which we call basic rules, and
we then tune the discharging process by introducing more complex rules.

\subsection{Basic rules for faces of size at least 12}\label{sec-large}

Each face $f_0$ of size $\ell\ge 12$ redistributes its charge as follows.

\begin{itemize}
\item Each A-triangle, B-triangle, and column with respect to $f_0$ receives \texttt{weak}$_{\ell}$ from $f_0$.
\item Each C-triangle and non-column $4$-face that shares an edge $v_1v_2$ with $f_0$
      receives $\text{\texttt{small}}_{\ell, a(v_1)}+\text{\texttt{small}}_{\ell, a(v_2)}$ from $f_0$,
      where $a(v_i)$ is $0$
      if the face sharing an edge with $f_0$ that incident with $v_i$ but not with $v_{3-i}$ is a $(\le\!4)$-face, and
      $a(v_i)$ is $1$, otherwise.
\item The sink of each isolated vertex incident with $f_0$ receives \texttt{iso}$_{\ell}$ from $f_0$ (recall that a sink may be a vertex or a face).
\end{itemize}

\noindent The amounts that are sent are defined in the following table.

\begin{center}
\begin{tabular}{|c|cccc|}
\hline
$\ell$ & \texttt{weak}$_{\ell}$ & \texttt{small}$_{\ell,0}$ & \texttt{small}$_{\ell,1}$ & \texttt{iso}$_{\ell}$\\
\hline 
$12$ & $\frac{4}{3}$ & $\frac{2}{3}$ & $\frac{1}{3}$ & $\frac{23827}{36960}$ \\
$13$ & $\frac{14023}{10080}$ & $\frac{14023}{20160}$ & $\frac{6137}{20160}$ & $\frac{1097}{1680}$\\
$(\ge\!14)$ & $2\bigl(1-\frac{4}{\ell}\bigr)$ & $1-\frac{4}{\ell}$ & $\frac{1}{2}\bigl(1-\frac{4}{\ell}\bigr)$ & $1-\frac{4}{\ell}$ \\
\hline
\end{tabular}
\end{center}

\subsection{Basic rules for faces of size between 5 and 11}\label{sec-med}

Fix an $\ell$-face $f_0$ with  $5\le \ell\le 11$.
The face $f_0$ sends $A_\ell$ to each incident A-triangle,
sends $B_\ell$ to each incident B-triangle, and sends $G_\ell$ to each incident column.
The amounts of charge that are sent are determined in the following table (note that $\ell\ge 6$
since a minimal counterexample does not contain TRIANGLE$_0$ or FOUR$_0$).

\begin{center}
\begin{tabular}{|c|ccc|}
\hline
$\ell$ & $A_\ell$ & $B_\ell$ & $G_\ell$ \\
\hline
6 & $\frac{1}{1}$ & $\frac{3}{4}$ & $\frac{3}{4}$ \\
7 & $\frac{1}{1}$ & $\frac{14}{15}$ & $\frac{9}{10}$ \\
8 & $\frac{1}{1}$ & $\frac{14}{15}$ & $\frac{82}{105}$ \\
9 & $\frac{17383}{15120}$ & $\frac{17383}{15120}$ & $\frac{2743}{2520}$ \\
10 & $\frac{8983}{7560}$ & $\frac{8983}{7560}$ & $\frac{16217}{15120}$ \\
11 & $\frac{4}{3}$ & $\frac{4}{3}$ & $\frac{4}{3}$ \\
\hline
\end{tabular}
\end{center}

Suppose that $uv_1v_2$ is a part of the boundary walk of $f_0$, and
the other face $f$ containing the edge $v_1v_2$ is either a C-triangle or a $4$-face.
Let $f'$ be the face incident with $uv_1$ distinct from $f_0$.
We define $t_f(v_1)$ as follows (when the face $f$ is clear from the context, we will omit the subscript).

\[ t_f(v_1)= 
  \begin{cases}
    0 & \mbox{if $\deg(v_1)\ge 4$ and $|{f'}|\geq 5$,} \\
    1 & \mbox{if $\deg(v_1)\ge 4$ and $|{f'}| \leq 4$, and} \\
    2 & \mbox{if $\deg(v_1)=3$ (and so $|{f'}| \geq 5$).}
  \end{cases}
\]

\noindent We next define $t(f)$ to be the value given by the following table.

\begin{center}
\begin{tabular}{|c|ccc|}
\hline 
\diagbox{$t(v_1)$}{$t(v_2)$}&0&1&2\\
\hline
0&0&1&2\\
1&1&3&4\\
2&2&4&5\\
\hline
\end{tabular}
\end{center}

If $f$ is a C-triangle,
then $f_0$ sends $C_{\ell,t(f)}$ to $f$ if $\ell\le 7$, and $C_{\ell,t(v_1)}+C_{\ell,t(v_2)}$ to $f$ if $\ell>7$.
Similarly, if $f$ is a non-column $4$-face,
then $f_0$ sends $D_{\ell,t(f)}$ to $f$ if $\ell\le 7$, and $D_{\ell,t(v_1)}+D_{\ell,t(v_2)}$ to $f$ if $\ell>7$.
The amounts of charge sent are given in the following table.

\begin{center}
\begin{tabular}{|c|cccccc|cccccc|}
\hline
$\ell$ & $C_{\ell,0}$ & $C_{\ell,1}$ & $C_{\ell,2}$ & $C_{\ell,3}$ & $C_{\ell,4}$ & $C_{\ell,5}$ & $D_{\ell,0}$ & $D_{\ell,1}$ & $D_{\ell,2}$ & $D_{\ell,3}$ & $D_{\ell,4}$ & $D_{\ell,5}$ \\
\hline
5 & $-\frac{11507}{36960}$ & $-\frac{7}{40}$ & $\frac{349}{840}$ & $-\frac{1}{7}$ & $\frac{13}{30}$ & $\frac{53}{120}$ & $\frac{1}{4}$ & $\frac{0}{1}$ & $\frac{349}{840}$ & $\frac{1}{15}$ & $\frac{1}{8}$ & $\frac{4}{7}$ \\
6 & $-\frac{10}{33}$ & $\frac{1}{336}$ & $\frac{1}{2}$ & $\frac{0}{1}$ & $\frac{1}{2}$ & $\frac{97}{160}$ & $\frac{1}{4}$ & $\frac{0}{1}$ & $\frac{3}{8}$ & $-\frac{13}{60}$ & $\frac{1}{8}$ & $\frac{67}{120}$ \\
7 & $\frac{2}{55}$ & $\frac{1}{336}$ & $\frac{211}{336}$ & $\frac{0}{1}$ & $\frac{1}{2}$ & $\frac{13}{15}$ & $\frac{1}{4}$ & $\frac{0}{1}$ & $\frac{3}{8}$ & $\frac{3}{7}$ & $\frac{3281}{20160}$ & $\frac{13}{15}$ \\
8 & $\frac{583}{1680}$ & $\frac{193}{840}$ & $\frac{7}{15}$ & & & & $\frac{41}{105}$ & $\frac{1}{4}$ & $\frac{1}{3}$ & & & \\
9 & $\frac{7223}{30240}$ & $\frac{1517}{7560}$ & $\frac{17383}{30240}$ & & & & $\frac{1009}{2160}$ & $\frac{5}{18}$ & $\frac{5017}{10080}$ & & & \\
10 & $\frac{83}{378}$ & $\frac{47851}{166320}$ & $\frac{8983}{15120}$ & & & & $\frac{20743}{40320}$ & $\frac{0}{1}$ & $\frac{4615}{8064}$ & & & \\
11 & $\frac{17}{33}$ & $\frac{7}{22}$ & $\frac{3}{5}$ & & & & $\frac{13}{22}$ & $\frac{26}{165}$ & $\frac{13}{22}$ & & & \\
\hline
\end{tabular}
\end{center}

Finally, suppose that $u_1vu_2$ is a part of the boundary of $f_0$ and $v$ is isolated.
For $i=1,2$, let $f_i$ be the face incident with $u_iv$ distinct from $f_0$.
If $|f_1|\ge r(\ell)$ and $|f_2|\ge r(\ell)$, then $f_0$ sends $E_{\ell,1}$ to the sink of $v$.
Otherwise, $f_0$ sends $E_{\ell,0}$ to the sink of $v$.
The values of $r(\ell)$, $E_{\ell,0}$ and $E_{\ell,1}$ are given in the following table.

\begin{center}
\begin{tabular}{|c|c|cc|}
\hline
$\ell$ & $r(\ell)$ & $E_{\ell,0}$ & $E_{\ell,1}$ \\
\hline
 5 & 15 & $\frac{1}{7}$ & $-\frac{61}{240}$ \\
 6 & 14 & $\frac{49}{240}$ & $-\frac{1}{15}$ \\
 7 & 13 & $\frac{79}{240}$ & $-\frac{1}{15}$ \\
 8 & 13 & $\frac{41}{105}$ & $\frac{9}{28}$ \\
 9 & 12 & $\frac{7}{15}$ & $\frac{1}{3}$ \\
10 & 12 & $\frac{7}{15}$ & $\frac{1}{3}$ \\
11 & 11 & $\frac{7}{15}$ & $\frac{1}{3}$ \\
\hline
\end{tabular}
\end{center}

\subsection{Basic rules for vertices of degree five and more}\label{sec-heavy}

In this subsection, we present basic rules for $(\ge\!5)$-vertices.
First, every $(\ge\!5)$-vertex sends $\frac{1}{4}$ to each incident $4$-face.

Every $5$-vertex incident with exactly one triangle $f$ and $m$ $4$-faces sends $\texttt{5\_to\_tri\_1}_m$ to $f$,
where $\texttt{5\_to\_tri\_1}_0=\frac{767}{1680}$, $\texttt{5\_to\_tri\_1}_1=\frac{737}{1680}$ and
$\texttt{5\_to\_tri\_1}_2=\frac{37}{120}$.
Finally, if a $5$-vertex $v$ is incident with faces $f_1,\ldots,f_5$ (in this order) and $|f_2|=|f_4|=3$,
then $v$ sends 
$\texttt{5\_to\_tri\_2\_light}=\frac{83}{140}$ to $f_2$ if $|f_1|\le 7$ and $|f_3|\le 7$, and
$\texttt{5\_to\_tri\_2\_heavy}=\frac{57}{140}$, otherwise.
This rule also applies to $f_4$ in the symmetric way.

Every $6$-vertex incident with exactly one triangle,
then it sends $\texttt{6\_to\_tri\_le2\_adj}=\frac{63}{80}$ to the incident triangle.
If a $6$-vertex is incident with two triangles, which cannot share an edge because of the reducibility of the configuration TFEDGE,
then it sends $\texttt{6\_to\_tri\_2\_opp}=\frac{767}{1680}$ to each of the two triangles.
Finally, if a $6$-vertex $v$ is incident with $f_1,\ldots,f_6$ (in this order) and $|f_2|=|f_4|=|f_6|=3$,
then $v$ sends $\texttt{6\_to\_tri\_3\_light}=\frac{113}{120}$ to $f_2$ if $\min(|f_1|,|f_3|)=5$ and $\max(|f_1|,|f_3|)\le 7$,
$v$ sends $\texttt{6\_to\_tri\_3\_all6}=\frac{8}{15}$ to $f_2$ if $|f_1|=|f_3|=6$, and
$v$ sends $\texttt{6\_to\_tri\_3\_heavy}=\frac{881}{1680}$ to $f_2$, otherwise.
The rule also applies to $f_4$ and $f_6$ in the symmetric way.
For example, if $|f_3|=|f_5|)=5$, then $v$ sends $\texttt{6\_to\_tri\_3\_all6}=\frac{8}{15}$ to $f_4$.

Finally, let $v$ be an $(\ge\!7)$-vertex and
let $f_1,\ldots,f_5$ be consecutive faces incident with $v$ such that $f_3$ is a $3$-face.
The vertex $v$ sends to $f_3$ $\texttt{6\_to\_tri\_3\_light}=\frac{113}{120}$ if $|f_1|=|f_5|=3$,
$\texttt{6\_to\_tri\_le2\_adj}=\frac{63}{80}$ if $\min\{|f_1|,|f_5|\}\le 4$ and $\max\{|f_1|,|f_5|\}\not=3$, and
$\texttt{6\_to\_tri\_2\_opp}=\frac{767}{1680}$, otherwise.

\subsection{Additional charge sent to 3-faces and 4-faces}\label{sec-addit}

Let $f_0$ be a $(\ge\!6)$-face and $u_1v_1v_2u_2$ be a part of its boundary walk.
Let $f_i$ be the other face incident with $u_iv_i$, $i=1,2$, and $f$ the other face incident with $v_1v_2$.
By symmetry, we can assume that $|f_1|\le |f_2|$.

If $f_0$ is a $6$-face, both $v_1$ and $v_2$ are $3$-vertices, $|f_i|\le \Delta^\star-1$ for $i=1,2$, and
$f$ is a non-column $4$-face, then $f_0$ sends $\texttt{light\_D\_extra} = \frac{1}{30}$ to $f$.

If $f_0$ is a $7$-face, both $v_1$ and $v_2$ are $3$-vertices, $|f_i|\le \Delta^\star-1$ for $i=1,2$, and
$f$ is a C-triangle, then $f_0$ sends $\texttt{light\_C\_extra}=\frac{1}{30}$ to $f$.

If $f_0$ is a $7$-face, $f$ is an A-triangle, and
$|f_1|=\Delta^\star-1$ (note that $|f_1|\ge\Delta^\star-1$ because of the absence of the configuration TRIANGLE$_0$
in a minimal counterexample),
then $f$ sends $\texttt{short\_to\_lightA}_{7,\Delta^\star-1,\Delta^\star-1}=\frac{1}{15}$ to $f$ if $|f_2|=\Delta^\star-1$, and
$f$ sends $\texttt{short\_to\_lightA}_{7,\Delta^\star-1,\Delta^\star}=\frac{1}{30}$ to $f$ if $|f_2|=\Delta^\star$.

If $f_0$ is a $8$-face, $f$ is an A-triangle, and $|f_1|\in\{\Delta^\star-2,\Delta^\star-1\}$,
then $f_0$ sends \texttt{short\_to\_lightA}$_{|f_1|,|f_2|}$ to $f$,
where the amounts are given by the following table.

\begin{center}
\begin{tabular}{|c|c|}
\hline
\texttt{short\_to\_lightA}$_{8,\Delta^\star-2,\Delta^\star-2}$ &$\frac{1}{7}$\\
\texttt{short\_to\_lightA}$_{8,\Delta^\star-2,\Delta^\star-1}$ &$\frac{1}{7}$\\
\texttt{short\_to\_lightA}$_{8,\Delta^\star-2,\Delta^\star}$ &$\frac{3}{28}$\\
\texttt{short\_to\_lightA}$_{8,\Delta^\star-1,\Delta^\star-1}$ &$\frac{1}{15}$\\
\texttt{short\_to\_lightA}$_{8,\Delta^\star-1,\Delta^\star}$ &$\frac{1}{30}$\\
\hline
\end{tabular}
\end{center}

If $f_0$ is a $9$-face, $f$ is an A-triangle, and $|f_1|=\Delta^\star-3$,
then $f_0$ sends $\texttt{face\_to\_lightA}_{9,2}=\frac{3257}{30240}$ to $f$ if $|f_2|=|f_1|$, and
$\texttt{face\_to\_lightA}_{9,1}=\frac{185}{6048}$, otherwise.

Finally, if $f_0$ is a $10$-face, $f$ is an A-triangle, and $|f_1|=\Delta^\star-4$,
then $f_0$ sends $\texttt{face\_to\_lightA}_{10,2}=\frac{583}{5040}$ to $f$ if $|f_2|=|f_1|$, and
$\texttt{face\_to\_lightA}_{10,1}=\frac{583}{10080}$, otherwise.

\subsection{Two-phase rules}\label{sec-th}

The rules described in this subsection have two phases:
first, some charge is sent to an edge $uv$ of $G$ and
then the edge $uv$ sends the received charge to one of the faces containg $uv$.
This way of describing discharging rules will be more convenient for the analysis of the sent charge in our proof.

Let $uv$ be an edge such that
both faces containing the edge $uv$ are $(\ge\!12)$-faces and $u$ is not a $3$-vertex contained in a $3$-face.
If $u$ is an $(\le\!4)$-vertex that is contained in exactly one $(\le\!4)$-face $f$,
then $f$ sends $\texttt{through\_heavy}=\frac{17}{80}$ to $e$.
If $u$ is a $4$-vertex contained in two $4$-faces or it is a $5$-vertex contained in two $3$-faces,
each of the two faces contaning $e$ sends $\frac{1}{2}\texttt{through\_heavy}=\frac{17}{160}$ to $e$.
Otherwise, $u$ sends $\texttt{through\_heavy}=\frac{17}{80}$ to $e$.
If $v$ is a $3$-vertex contained in a $3$-face $vv'v''$ and
the other face $f'$ containing the edge $v'v''$ is a $(\le\!11)$-face,
then $e$ sends $\texttt{through\_heavy}=\frac{17}{80}$ to $f$.

Let $uv$ be an edge such that both faces containing $uv$ has size between $5$ and $10$ (inclusively).
Let $f$ be one of the faces containing $uv$ and $u'$ the neighbor of $u$ incident with $f$ that is different from $v$.
If the edge $uu'$ is contained in a $3$-face $f'$, and
either $u$ is an $(\ge\!6)$-vertex or $u$ is a $5$-vertex contained in only one $3$-face,
then $u$ sends $\texttt{through\_heavy}=\frac{17}{80}$ to $e$, which then sends $\texttt{through\_heavy}=\frac{17}{80}$ to $f'$.
Note that if $u$ is a $(\ge\!6)$-vertex, then this rule may apply twice,
triggered by each of the two faces containing the edge $uv$.

Finally, let $uv$ be an edge contained in a face of size between $5$ and $10$ (inclusively) and a $(\ge\!12)$-face $f$.
Let $u'$ be the neighbor of $u$ incident with $f$ that is different from $v$.
If the edge $uu'$ is contained in a $3$-face $f'$ and $u$ is a $5$-vertex contained in two $3$-faces,
then $f'$ sends $\texttt{through\_heavy}=\frac{17}{80}$ to $e$,
which then sends $\texttt{through\_heavy}=\frac{17}{80}$ to the $3$-face containing $u$ that is different from $f'$.

\subsection{Additional special rules}\label{sec-spec}

Let a $4$-face $f=v_1v_2v_3v_4$ and a $5$-face $f'$ share the edge $v_1v_2$.
If $v_3$ is an $(\ge\!4)$-vertex, $v_4$ is an $(\ge\!4)$-vertex or $v_3v_4$ is contained in a $(\ge\!6)$-face,
then $f$ sends $\texttt{four\_to\_five}=\frac{109}{840}$ to $f'$.

If $f$ and $f'=v_1v_2v_3v_4$ are two $4$-faces sharing the edge $v_1v_2$,
both $v_1$ and $v_2$ are $4$-vertices, and the other faces containing $v_1v_4$ and $v_2v_3$ are also $4$-faces,
then $f$ sends $\texttt{four1}=\frac{1}{2}$ to $f'$.

If $f=v_1v_2v_3v_4$ is a $4$-face, $v_1$ is a $4$-vertex contained in a $3$-face $f'$, and
the other faces containing $v_1v_4$ and $v_2v_3$ are $(\ge\!\Delta^\star-1)$-faces,
then $f$ sends $\texttt{four2}=\frac{1}{2}$ to $f'$.

If $v_1v_2v_3$ is a part of the boundary walk of an $\ell$-face $f$, $\ell\in\{5,11\}$,
$v_1$ is a $3$-vertex, and
both $v_1v_2$ and $v_2v_3$ are incident with C-triangles,
then the $3$-face containing $v_1v_2$ sends $\texttt{$\star$\_CC\_to\_5\_extra}=\frac{37}{240}$ to $f$ if $f$ is a $5$-face, and
$\texttt{$\star$\_CC\_to\_11\_extra}=\frac{14}{165}$, otherwise (when $f$ is a $11$-face).

If $v_1v_2v_3$ is a part of the boundary walk of a 10-face $f$,
$v_2v_3$ is contained in an A-triangle $f'$, and
the other face containing $v_1v_2$ is a $(\ge\!13)$-face,
then $f$ sends $\texttt{10\_to\_13\_A\_extra}=\frac{89}{6048}$ to $f'$.

Finally, if $v_1v_2v_3$ is a part of the boundary walk of an 11-face $f$,
both $v_1v_2$ and $v_2v_3$ are contained in faces of size $5$ or $6$,
$v_2$ is a $4$-vertex contained in a $3$-face $f'$,
then $f$ sends $\texttt{11\_to\_opp\_66tri\_extra}=\frac{28}{165}$ to $f'$,
which is the sink of $v_2$.

\section{Analysis of final charges}
\label{sec-analysis}

In this section,
we argue that if a graph $G$ that does not contain any of the reducible configurations (which were identified in Section~\ref{sec-redu}),
is assigned charge as described at the beginning of Section~\ref{sec-rules} and
then this charge is redistributed using the rules described in the rest of Section~\ref{sec-rules},
then the final charge of each vertex, edge and face of $G$ is non-negative.
Since the charge is preserved by the rules and the initial amount of charge was negative,
this contradicts the existence of a counterexample to Theorem~\ref{main-thm}.

The final charge of edges is easy to analyze.
The edges are only involved in the rules described in Subsections~\ref{sec-th} and
each edge sends out as much as it has received.

The analysis of the final amount of charge of vertices and faces is more involved.
We performed the analysis with the assistance of a computer.
The program is available at \url{http://www.ucw.cz/~kral/cyclic-16/} as
the file \texttt{test-discharging.lhs}.
We used Literate Haskell to prepare the program:
compiling the file with Latex produces a detailed description of how the program works, and
compiling it with GHC produces an executable file that performs the analysis.
The former file is available on the webpage, too.

In the rest of the section, the rules are referred to by the names of the constants described the amount of charge transferred.
For example, the \texttt{iso} rules are the rules described in the third point in Section~\ref{sec-large}.

\subsection{Final charge of vertices}

We now give details how the amount of final charge of vertices is analyzed.
Since $G$ is $3$-connected, its minimum degree is at least three.
If a $3$-vertex $v$ is contained in a $(\le\!4)$-face,
then it gets $1$ unit of charge from the incident $(\le\!4)$-face(s) and is not affected by any other rules.
If a $3$-vertex $v$ is not contained a $(\le\!4)$-face,
then it receives charge described by \texttt{iso} and \texttt{E} rules from Sections~\ref{sec-large} and \ref{sec-med}, and
it can send out charge by the \texttt{through\_heavy} rules from Section~\ref{sec-th}.
In particular, the amounts received and sent only depend on the sizes of the faces containing $v$.
Hence, the program just enumerates all possibilities and checks that the final charge of $v$ is non-negative.
We proceed similarly for $4$-vertices, $5$-vertices, $6$-vertices and $7$-vertices.
Note that a $4$-vertex contained in a $(\le\!4)$-face is unaffected by any rules (its sink is the incident $(\le\!4)$-face,
so it does not receive any charge by the \texttt{iso} and \texttt{E} rules),
so such vertices need not be analyzed.

Consider a $d$-vertex $u$, $d\ge 8$, and let $f_1,\ldots,f_d$ be the faces incident with $u$ (in this cyclic order).
For $1\le i\le d$, define $c_i$ to be the following charge.
If $|f_i|=3$, then $c_i$ is
the amount of charge sent from $u$ to $f_i$ by the rules from Section~\ref{sec-heavy} plus
the amount of charge sent from $u$ to an edge $uv$ by the \texttt{through\_heavy} rules described in the last two paragraphs of Section~\ref{sec-th}.
If $|f_i|=4$, $c_i$ is the amount of charge $u$ sends to $f_i$ by the rules from Section~\ref{sec-heavy}.
Otherwise, $c_i$ is half the amount of charge sent by the \texttt{through\_heavy} rules described in the second paragraph of Section~\ref{sec-th} minus
the amount of charge received from $f_i$ by the \texttt{iso} and \texttt{E} rules.
Observe that $c_i$ depends only on the sizes of the faces $f_{i-2},\ldots,f_{i+2}$ (with indices modulo $d$) and
let $q_i=\frac{1}{2}c_{i-1}+c_i+\frac{1}{2}c_{i+1}$ (again with indices modulo $d$).
The program enumerates all possible sizes of the faces $f_{i-3},\ldots,f_{i+3}$ and checks that $q_i\le 1$.
This yields that
$$\sum_{i=1}^d c_i=\frac{1}{2}\sum_{i=1}^d q_i\le \frac{d}{2}\;\mbox{.}$$
Hence, the total amount of charge sent out by $u$ is at most $d/2\le d-4$ and
its final charge is non-negative.

\subsection{Final charge of faces}

The final amounts of charge of faces are analyzed in a different way depending on the face sizes.
Let us start by considering a $3$-face $f=v_1v_2v_3$, and
let $f_i$ be the other face containing the edge $v_iv_{i+1}$ (indices modulo three).
The \emph{shape} of a face $f$ consists of the information on the sizes of the faces $f_1$, $f_2$ and $f_3$ and
the information whether the faces cyclically adjacent to $f_i$ at $v_i$ and $v_{i+1}$ are $3$-faces, $4$-faces or $(\ge 5)$-faces.
The shape fully determines the amount of charge sent by $f$ to incident $3$-vertices and
the amount of charge received from the incident faces by the basic rules from Subsections~\ref{sec-large} and~\ref{sec-med}, and
the amount of charge received by the rules from Subsection~\ref{sec-addit} and~\ref{sec-spec} except for the rule\texttt{11\_to\_opp\_66tri\_extra}.
Let $c_0$ be the total amount of this charge.

The charge not accounted for in $c_0$ is sent by the \texttt{E} rules through $4$-vertices,
by the rules from Subsections~\ref{sec-heavy} and~\ref{sec-th} and the rule \texttt{11\_to\_opp\_66tri\_extra}.
Each of these rules can be associated with one of the vertices $v_i$, $i=1,2,3$, and
its amount only depends on the sizes of the faces containing $v_i$ (in addition to the shape of $f$).
Hence, we can determine the worst case charge $c_i$ for each vertex $v_i$ independently of the other two vertices of $f$.
We then verify that $c_0+c_1+c_2+c_3-1$ is non-negative for each possible shape of a $3$-face.

The analysis of the final charge of $4$-faces is similar to that of $3$-faces.
We now focus on faces with sizes between $5$ and $13$ (inclusively).
The \emph{inventory} of a face is the number of adjacent A-triangles, B-triangles, columns and $(\le\!4)$-faces
distinguished by the number of their vertices that are incident with another $(\le\!4)$-face.
The inventory is enough to determine the final charge of $(\ge\!12)$-faces.
The program enumerates all possible inventories of $(\ge\!12)$-faces and
checks that the final charge of all $(\ge\!12)$-faces is non-negative.

The program also enumerate possible inventories of $\ell$-faces, $5\le\ell\le 11$,
discards those that give a non-negative lower bound on the final amount of charge of the considered face $f$.
For each of the non-discarded inventories,
the program enumerates all cyclic orders determining which edges of the $\ell$-face are contained in the elements of the inventory.
Some of the enumerated configurations can be excluded by the reducible configurations and get discarded (this actually finishes off
the analysis of $6$-faces and $11$-faces).
In addition, lower bounds on the sizes of the other faces adjacent to $f$ are obtained,
e.g., the configuration GEN$_2$ is used to establish that an incident face must have size at least $\Delta^\star+7-\ell$.
In the remaining cases, 
the program enumerates all possible sizes of faces that affect the charge sent or received by $f$,
i.e., faces next to A-triangles, faces incident with $3$-vertices of C-triangles at $(\le\!7)$-faces, and
faces incident with the $3$-vertices of non-column $4$-faces at $5$-faces, and
it checks that the final charge of $f$ is non-negative.

It remains to analyze the final amount of charge of $(\ge\!14)$-faces.
We account the charge sent out by an $\ell$-face $f$ to its vertices, $\ell\ge 14$.
If $v_1v_2$ is an edge of $f$ contained in an A-triangle, a B-triangle or a column,
then $\texttt{weak}_{\ell}/2=1-\frac{4}{\ell}$ is assigned to each of $v_1$ and $v_2$.
If $v_1$ is an isolated vertex, it is assigned \texttt{iso}$_{\ell}=1-\frac{4}{\ell}$.
If $v_1v_2v_3$ is a path on the boundary of $f$ and $v_1v_2$ is contained in a C-triangle or a non-column $4$-face,
then $\texttt{small}_{\ell, a(v_1)}$ is assigned to $v_1$ and $\texttt{small}_{\ell, a(v_2)}$ is assigned to $v_2$.
If the edge $v_2v_3$ is also in a triangle or a $4$-face $f$ (which cannot be an A-triangle, B-triangle or a column),
then $a(v_2)=1$ and we assign $\texttt{small}_{\ell,1}$ to $v_2$ in addition,
i.e., $v_2$ is assigned $2\texttt{small}_{\ell, 1}=1-\frac{4}{\ell}$ in total.
Otherwise, $a(v_2)=0$ and $v_2$ is assigned $\texttt{small}_{\ell, 0}=1-\frac{4}{\ell}$.
We conclude that the charge sent out by $f$ is at most $\ell\left(1-\frac{4}{\ell}\right)=\ell-4$,
i.e., the final charge of $f$ is non-negative.

\section*{Acknowledgements}

The authors would like to thank Jan van den Heuvel for pointing out a flaw in one of their arguments.
They would also like to thank the two anonymous reviewers for their careful reading of the manuscript and their comments,
which helped to improve the quality of the presentation susbtantially.

\end{document}